\documentclass[a4paper, 11pt]{article}
\usepackage{bbm}
\usepackage{amsfonts}
\usepackage{amsmath}
\usepackage{amsthm,mathrsfs}
\usepackage{amsfonts,amssymb}
\usepackage{float}
\usepackage{pdfsync}
\usepackage{graphicx}
\usepackage{caption}
\usepackage{hyperref}
\allowdisplaybreaks[4]

\newtheorem{theorem}{Theorem}[section]
\newtheorem{proposition}{Proposition}[section]
\newtheorem{lemma}{Lemma}[section]

\newtheorem{remark}{Remark}[section]
\theoremstyle{remark}
\theoremstyle{corollary}

\theoremstyle{remark}
\theoremstyle{step}

\numberwithin{equation}{section}

\def\dfrac{\displaystyle\frac}

\def\iy{\infty}

\def\p{\partial}
\def\s{\sigma}
\def\d{\delta}
\def\v{\varphi}
\def\k{\kappa}
\def\l{\lambda}
\def\L{\Lambda}
\def\g{\gamma}
\def\O{\Omega}
\def\t{\theta}
\def\ss{\setminus}

\def\rr{{\mathbb R}}

\begin{document}
\title{{$H^{1+\alpha}$ estimates for the fully nonlinear parabolic thin obstacle problem}
\footnote
{Supported by the National Natural Science Foundation of China (No.11771023).}}

\date{}\maketitle
\vspace{-2cm}
\begin{center}
Xi Hu and Lin Tang

\end{center}
{\bf Abstract:}  We study the regularity of the viscosity solution to the fully nonlinear parabolic thin obstacle problem. In particular, we prove that the solution is local $H^{1+\alpha}$ on each side of the smooth obstacle, for some small $\alpha>0.$ Following the method which was first introduced for the harmonic case by Caffarelli in 1979, we extend the results of Fern\'{a}ndez-Real (2016) who treated the fully nonlinear elliptic case. Our results also extend those of Chatzigeorgiou (2019) in two ways. First, we do not assume solutions nor operators to be symmetric. Second, our estimates are local, in the sense that do not rely on the boundary data.
\vspace{0.2cm}\\
{\bf Keywords:} Parabolic thin obstacle problem; Fully nonlinear parabolic equations; Free boundary problems
\vspace{0.2cm}\\
{\bf Mathematics Subject Classification:} 35R35, 35R45, 35K55

\section{Introduction}
The thin obstacle problem originated in the context of elasticity. We can think of the thin obstacle problem as an equation describing the shape of an elastic membrane that is pushed from below by a very thin object (of co-dimension 1). Important cases of obstacle type problems occur when the operators involved are fractional powers of the Laplacian as well as nonlinear operators, since they appear in the study of the quasi-geostrophic flow model \cite{CP}, anomalous diffusion \cite{BG} and American options with jump processes \cite{Sl}.

The study of the regularity properties of the thin obstacle problem has a long history, dating back to the works of Lewy \cite{Le} and Richardson \cite{Rd}, who established the optimal regularity result for the two-dimensional problem by means of complex analysis techniques. Subsequently, the full higher dimensional problem was studied in different generalities by various authors; see \cite{Fj,C1,Kd,UN}.
In particular, the $C^{1,\alpha}$ regularity of the weak solutions for the harmonic case and the divergence case for regular enough coefficients was proved in \cite{C1}.  Results for more general divergence form elliptic operators can be found in \cite{UN}. Since then, the thin obstacle problem was extensively studied in the elliptic case. Recently, in their seminal articles, Athansopoulos and Caffarelli \cite{AC} and Athansopoulos, Caffarelli and Salsa \cite{ACS} introduced Almgren's frequency function as a powerful tool of obtaining optimal regularity estimates of solutions and  the $C^{1,\alpha}(\alpha\in (0,1])$ regularity of the free boundary for the Laplacian equation. Later this was extended by Silvestre in \cite{Sl} and Caffarelli, Salsa and Silvestre in \cite{CSS} to the related obstacle problem for the fractional Laplacian. While the constant coefficient situation was studied thoroughly subsequently \cite {GN, GP1, PS, Fo}, the variable coefficient situation with low coefficient regularity was treated in \cite {GS, GP, KR1, KR}. More recently, for fully nonlinear elliptic operators, the regularity of the viscosity solution was first proved by Milakis and Silvestre in \cite{MS}. Such a result had been generalized to the general nonsymmetric setting by Fern\'{a}ndez-Real  in \cite{FX}. An important difficulty in the study of the free boundary problem is lack of monotonicity formulas for fully nonlinear operators. The regularity of the free boundary for fully nonlinear operators was solved by Ros-Oton and Serra in \cite{RX}. Their proofs are completely independent from those in \cite{ACS}, and do not use any monotonicity formula.

The corresponding regularity theory for thin obstacle problems of parabolic type is much less developed. The $H^{1+\alpha}$ regularity of the weak solution for the case of heat equation and the case of smooth enough linear parabolic equation was obtained  in 1982 by Athanasopoulos in \cite{A}. The case of more general linear parabolic operators was examined by Uraltseva in \cite{UN1} and Arkhipova in \cite{A1}. Recently, Danielli, Garofalo, Petrosyan and To \cite{DG} generalize Almgren's frequency formula to the parabolic case and obtained the optimal and free boundary regularity for heat equation. For the case of parabolic fractional Laplacians and degenerate parabolic equations we refer the reader to \cite {CF, ACM, BD, BD1}. Only very recently, Chatzigeorgiou \cite{CG} obtained the regularity of the viscosity solution to the thin obstacle problem for fully nonlinear parabolic operators, which extended the results of \cite{MS} to the parabolic case.

A natural generalization of the results of \cite{CG} is to consider the general nonsymmetric setting. This is the problem we study in this paper. On one hand, our results extend the results of \cite{FX} who treated the fully nonlinear elliptic case.  On the other hand, our results also extend the results of \cite{CG} in two ways. First, we do not assume anything on the boundary data, so that we give a local estimate. Second, we also consider nonsymmetric solutions $u$ with operators not necessarily satisfying any symmetry assumption, and prove $H^{1+\alpha}$ regularity for such solutions. In future work, we plan to extend the results of Ros-Oton and Serra (2017) who treated free boundary regularity for the fully nonlinear elliptic operators to the parabolic case.


Given a domain $D$ in $\rr^{n+1}$ with an $n$-dimensional surface $S$ inside, the parabolic thin obstacle problem involves a function $u: D\rightarrow \rr,$ an obstacle $\v: S\rightarrow \rr$, a Dirichlet boundary condition given by $u_0: \p_pD\rightarrow\rr$, and a second order elliptic operator $L$,
\begin{equation*}\label{1}
    \left\{
   \begin{array}{ll}
Lu-\p_tu=0\quad  &\text{in}~D\ss\{(x,t)\in S: u(x,t)=\v(x,t)\}\\
Lu-\p_tu\leq 0&\text{in}~D\tag{1.1}\\
u\geq \v&\text{on}~S \\
u=u_0 &\text{on}~\p_pD.
 \end{array}
 \right.
\end{equation*}

Here, we study a nonlinear version of problem (1.1). More precisely, we study (1.1) with $Lu=F(D^2u)$, a convex fully nonlinear uniformly elliptic operator. Since all of our estimates are of local character, we consider the problem in $Q_1,$
\begin{equation*}\label{1}
    \left\{
   \begin{array}{ll}
F(D^2u)-\p_tu=0\quad  &\text{in}~Q_1\setminus \{u=\v\}\\
F(D^2u)-\p_tu\leq 0&\text{in}~Q_1\tag{1.2}\\
u\geq \v&\text{on}~Q_1\cap\{x_n=0\}.\\
 \end{array}
 \right.
\end{equation*}
Here, $\v: Q_1\cap\{x_n=0\}\rightarrow \rr$ is the obstacle, and we assume that it is $H^2.$ We also assume that
\begin{align*}
&F~\text{is convex, uniformly elliptic}\\
&\text{with ellipticity constants}~ 0<\l\leq \L,~ \text{and with}~ F(0)=0\tag{1.3}.
\end{align*}

When $u$ is symmetric, the following problem was studied by Chatzigeorgiou in \cite{CG}
\begin{equation*}\label{1}
    \left\{
   \begin{array}{ll}
F(D^2u)-\p_tu=0\quad \quad &\text{in}~Q^+_1\\
\max\{u_{x_n},\v-u\}=0&\text{on}~Q_1\cap\{x_n=0\}\tag{1.4}\\
u=u_0 &\text{on}~\p_pQ^+_1\setminus \{{Q_1\cap\{x_n=0\}\}},\\
 \end{array}
 \right.
\end{equation*}
where $u_t$ is locally bounded by above in $Q^+_1$ in subsection 2.2. In addition, they also implicitly assume a symmetry condition on the operator $F$. This assumption is important in \cite{CG} to prove semiconvexity of solutions.

In the linear case, one can symmetrise solutions to (1.2), and then the study of such solutions reduces to problem (1.4). However, in the present nonlinear situation an estimate for (1.4) does not imply one for (1.2).

Now we state our main result as follows.
\begin{theorem}
Let $F$ be a fully nonlinear operator satisfying $(1.3)$. Let $u$ be any viscosity solution to $(1.2)$ with $\v\in H^2$ and $\p_tu$ is locally upper bounded. Then, $u\in H^{1+\alpha}(\overline{Q^+_{1/2}})\cap H^{1+\alpha}(\overline{Q^-_{1/2}})$ and, $$\|u\|_{H^{1+\alpha}(\overline{Q^+_{1/2}})}+\|u\|_{H^{1+\alpha}(\overline{Q^-_{1/2}})}\leq C\big(\|u\|_{L^\infty(Q_1)}+\|\v\|_{H^2(Q_1\cap\{x_n=0\})}\big)$$
for some constants $\alpha>0$ and $C$ depending only on $n,\l$ and $\L.$
\end{theorem}
\begin{remark}
In particular, for the Laplacian or fractional Laplacians, under some smoothness assumption on initial data $u(x,0)=u_0(x)$, the boundedness of $u_t$ is obtained by using a penalization method; see \cite{A} and \cite{CF}.
\end{remark}

The paper is organized as follows. In Section 2 we give a list of notations used throughout this paper. In Section 3 we examine Lipschitz continuity in the space variables. In addition, we prove the semiconvexity of solutions. Our proof of the semiconvexity of solutions is completely different from the one done in \cite{CG} and follows by means of a Bernstein's technique. Finally in Section 4 we prove Theorem 1.1. More precisely, to prove Theorem 1.1, we first define a symmetrised solution to the problem and follow the steps in \cite{C1} and \cite{FX} using approximate inequalities satisfied by the symmetrised solution, which yields the regularity of the symmetrised normal derivative at free boundary points. Then, applying the ideas from \cite{W2}, we show that the $H^{1+\alpha}$ regularity of the original function $u$ at free boundary points. Finally, we show that the regularity of $u$ at free boundary points yields the regularity of the symmetrised normal derivative at all points on $\{x_n=0\}$, and this yields the regularity of $u$ on either side of the obstacle.

\section{Notations}
The Euclidean ball in $\rr^n$ and the elementary cylinder in $\rr^{n+1}$ will be denoted by $$B_r(x_0):=\{x\in \rr^n:|x-x_0|<r\},\quad Q_r(x_0,t_0):=B_r(x_0)\times(t_0-r^2,t_0],$$
respectively. We define the following half and thin balls in $\rr^n$, for $r>0$, $x_0\in \rr^{n-1},$
$$B^{\pm}_r(x_0):=B_r(x_0)\cap\rr^n_{\pm},\quad B^*_r(x_0):=B_r(x_0)\cap\rr^{n-1}$$ and the following half and thin cylinders in $\rr^{n+1}$, for $r>0,$ $x_0\in \rr^{n-1}$, $t_0\in \rr,$
$$Q^{\pm}_r(x_0,t_0)=B^{\pm}_r(x_0)\times(t_0-r^2,t_0],\quad Q^*_r=B^*_r(x_0)\times(t_0-r^2,t_0].$$
Moreover, we call the coincidence set $$\triangle^*=\{(x',t)\in Q^*_1:u(x',0,t)=\v(x',t)\},\quad \triangle=\triangle^*\times\{0\},$$
and its complement in $Q^*_1$ is denoted by $$\O^*=Q^*_1\setminus \triangle^*,\quad \O=\O^*\times\{0\}.$$
Notice that $\O^\circ$, $\overline{\O}$ and  $\p \O$ will be the interior, the closure and the boundary of the domain $\O\subset \rr^{n+1},$ respectively, in the sense of the Euclidean topology of $\rr^{n+1}.$ Let us define the parabolic interior to be, $$\mathrm{int}_p(\O):=\{(x,t)\in \rr^{n+1}: \text{there exists }~ r>0~ \text{so that } Q^\circ_r(x,t)\subset\O\}$$
and the parabolic boundary, $\p_p(\O):=\overline{\O}\setminus \mathrm{int}_p(\O).$ In addition, the parabolic distance will be denoted by $$ p(P_1,P_2):=\max\{|x-y|,|t-s|^{1/2}\} \quad\text{for}~P_1(x,t), P_2(y,s)\in \rr^{n+1}.$$

Next, we define the corresponding parabolic H\"{o}lder spaces as follows. Given $\beta\in (0,2]$ and $f$ defined on $\Omega\subset \mathbb{R}^{n+1}$, we denote $$\langle f\rangle_{\beta;(x_0,t_0)}:=\sup\bigg\{\dfrac{|f(x_0,t)-f(x_0,t_0)|}{|t-t_0|^{\beta/2}}:(x_0,t)\in \Omega \setminus \{(x_0,t_0)\}\bigg\}\quad \text{for}~(x_0,t_0)\in \Omega.$$
And we denote $\langle f\rangle_{\beta;\Omega}:=\displaystyle\sup_{(x_0,t_0)\in \Omega }\langle f\rangle_{\beta;(x_0,t_0)}.$ For any $l>0$, we write  $l=k+\alpha,$ where $k$ is a nonnegative integer and $\alpha\in (0,1],$ and we define
$$\langle f\rangle_{l;\Omega}:=\sum_{|\beta|+2j=k-1}\langle D^\beta_xD^j_t f\rangle_{\alpha+1},$$

$$[f]_{l;\Omega}:=\sum_{|\beta|+2j=k}[ D^\beta_xD^j_t f]_{\alpha},$$

$$|f|_{l;\Omega}:=\sum_{|\beta|+2j\leq k}|D^\beta_xD^j_t f|+[f]_l+\langle f\rangle_l,$$
and we let $H^l(\Omega):=\{f:|f|_l<\infty\}$.

Hereafter, when we say that a constant is universal, we mean that it depends only on any or all of $n,\l,\L$ or $\alpha.$ Throughout this paper, we will denote $K:=\|u\|_{L^\infty(Q_1)}+\|\v\|_{H^2(Q^*_1)}$. Unless otherwise stated, $c~\text{and}~C$ will also denote positive universal constants that may change from line to line.
\section{Lipschitz estimate and semiconvexity properties}
\subsection{Lipschitz estimate}
We begin with a proposition showing any solution to (1.2) is spatial Lipschitz.
\begin{proposition}
Let $u$ be any solution to $(1.2)$ with $F$ satisfying $(1.3)$ and $\varphi\in H^{2}$. Then $u$ is spatial Lipschitz in $Q_{1/2}$ with,
$$\|u\|_{\mathrm{Lip}(Q_{1/2})}\leq C\Big(\|u\|_{L^\infty(Q_1)}+\|\v\|_{H^{2}(Q^*_1)}\Big),\eqno(3.1)$$
for some universal constant $C$.
\end{proposition}
\begin{proof}
We thicken the obstacle $\v$. We extend the obstacle $\v$ to a function $h$ defined in the whole $Q_1,$ and we treat $u$ as a solution to a classical ``thick'' obstacle problem. Namely, we consider the viscosity solutions to
\begin{equation*}\label{1}
    \left\{
   \begin{array}{ll}
F(D^2h)-\p_th=0\qquad &\text{in}~Q^+_1\\
h(x',0,t)=\v(x',t)& \text{on}~ Q^*_1\tag{3.2}\\
h=-\|u\|_{L^\infty(Q^+_1)}&\text{on}~\p_pQ^+_1\setminus Q^*_1,
 \end{array}
 \right.
\end{equation*}
and analogously
\begin{equation*}\label{1}
    \left\{
   \begin{array}{ll}
F(D^2h)-\p_th=0\qquad &\text{in}~Q^-_1\\
h(x',0,t)=\v(x',t)&\text{on}~ Q^*_1\tag{3.3}\\
h=-\|u\|_{L^\infty(Q^-_1)}& \text{on}~\p_pQ^-_1\setminus Q^*_1.
 \end{array}
 \right.
\end{equation*}

Note that $h$ is Lipschitz in $Q_{7/8}$; see \cite{CG1} and \cite{LG}. Then
we have $$\|h\|_{\mathrm{Lip}(Q_{7/8})}\leq CK.$$
Moreover, using  maximum principle we can obtain that  $u\geq h$ in $Q_1.$ In addition, $u$ is a solution to a classical obstacle problem in $Q_1$ with $h$ as the obstacle.
We show next that this implies $u$ is spatial Lipschitz, with a quantitative estimate.

Since $h$ is Lipschitz, fixed any $(x_0,t_0)\in Q_{1/2}$ and $0<r<1/4$, there exists universal constant $C_0$ such that $$\sup_{Q_r(x_0,t_0)}|h(x,t)-h(x_0,t_0)|\leq C_0Kr.\eqno(3.4)$$

By the strong maximum principle, the coincidence set $\{u=h\}$ is $\triangle$, the coincidence set of parabolic thin obstacle problem. Assume that $(x_0,t_0)\in \triangle$. Since $u\geq h$, $$\inf_{Q_r(x_0,t_0)}(u(x,t)-u(x_0,t_0))\geq -C_0Kr.\eqno(3.5)$$

Let $q(x,t)=u(x,t)-u(x_0,t_0)+C_0Kr.$ Using (3.5), we know $q\geq 0$ in $Q_r(x_0,t_0).$ In addition, by (3.4), $$q(x,t)\leq 2C_0Kr\quad\text{on}~Q_r(x_0,t_0)\cap\triangle.$$ It is easy to check that $q$ is a supersolution of the equation $$F(D^2q)-\p_tq=F(D^2u)-\p_tu\leq 0\quad\text{in}~Q_r(x_0,t_0).$$

Let $\bar{q}$ be the viscosity solution to $F(D^2\bar{q})-\p_t\bar{q}=0$ in $Q_r(x_0,t_0)$ with $\bar{q}=q$ on $\p_pQ_r(x_0,t_0).$ We have $\bar{q}\leq q$ in $Q_r(x_0,t_0)$. Applying the nonnegativity of $\bar{q}$ on the boundary, $\bar{q}\geq 0$ in $Q_r(x_0,t_0).$

Therefore, $q<\bar{q}+2C_0Kr$ on $\p_p Q_r(x_0,t_0),$ and $q\leq \bar{q}+2C_0Kr$ in $Q_r(x_0,t_0)\cap \triangle.$ Thus, $$q\leq \bar{q}+2C_0Kr\quad \text{in}~Q_r(x_0,t_0).$$

On the other hand, we know that $0\leq \bar{q}(x_0,t_0)\leq q(x_0,t_0)=C_0Kr.$ Applying Harnack inequality, we obtain $\bar{q}\leq CC_0Kr$ in $Q_{r/2}(x_0,t_0).$ Furthermore, we deduce $$u(x,t)-u(x_0,t_0)\leq CC_0Kr,\eqno(3.6)$$ for some constant $C>0.$ Therefore, combining (3.5) with (3.6), $$\sup_{Q_r(x_0,t_0)}|u(x,t)-u(x_0,t_0)|\leq C\Big(\|u\|_{L^\infty(Q_1)}+\|\v\|_{H^{2}(Q^*_1)}\Big)r,$$ for some universal constant $C$.

For the points of the coincidence set, we have obtained that the solution is spatial Lipschitz. Applying interior estimates, we show next that the solution is spatial Lipschitz inside $Q_{1/2}$.

Take any points $P_1=(x,t), P_2=(y,s)\in Q_{1/2},$ and let $r=p(P_1,P_2).$ Define $$\rho:=\min\{\text{dist}_p(P_1,\triangle),\text{dist}_p(P_2,\triangle)\},$$
and let $P^*_1,P^*_2\in \triangle, P^*_1=(x',0,\bar{t}), P^*_2=(y',0,\bar{s})$ for $(x',\bar{t}), (y',\bar{s})\in \triangle^*$, be such that $\text{dist}_p(p_1,\triangle)=p(P_1,P^*_1)$ and $\text{dist}_p(P_2,\triangle)=p(P_2,P^*_2).$ We split into two cases.\\
\\
$\mathrm{Case}~1:$ If $\rho\leq 4r,$ then
\begin{align*}
|u(P_1)-u(P_2)|&\leq |u(P_1)-u(P^*_1)|+|u(P_2)-u(P^*_2)|+|\v(x',\bar{t})-\v(y',\bar{s})|\\
&\leq C\rho+C(r+\rho)+2C(r+\rho)\\
&\leq Cr
\end{align*}
for some constant $C.$ We are using here that $\v\in H^2$ and that if $p(P_1,P^*_1)=\rho,$ then $p(P_2,P^*_2)\leq r+\rho$ and $p(P^*_1,P^*_2)\leq 2(r+\rho).$\\
\\
$\mathrm{Case}~2:$ If $\rho>4r,$ we can use interior estimate. Assume $(x,t)$ is such that $\text{dist}_p((x,t),\triangle)=\rho,$ and note that $Q_{\rho/2}(x,t)\subset Q_1\setminus\triangle$, so that $F(D^2u)-\p_tu=0~\text{in}~Q_{\rho/2}(x,t)$. Now, we can use the interior Lipschitz estimates (see \cite{W1}) and the supremum and the infimum of $u$ in $Q_{\rho/2}(x,t)$ controlled by $C\rho+\v(P^*_1)$ and $-C\rho+\v(P^*_1)$,$$[u]_{\text{Lip}(Q_{\rho/4})}\leq \frac{C}{\rho}{\mathop{\text{osc}}\limits_{Q_{\rho/2}(x,t)}}u\leq C$$ for some constant $C.$ Thus, we finish the proof of Proposition 3.1.
\end{proof}
\subsection{Semiconvexity and semiconcavity estimates}
Before continuing to prove the semiconvexity and semiconcavity results, we adapt a change of variables that will be useful. Following \cite{FX}, throughout the paper we will assume that there exists a fixed symmetric uniformly  elliptic operator $\hat{L}$ such that $$\hat{L}^{ij}\p_{x_ix_j}w\leq F(D^2w),\quad \text{and}~\hat{L}^{in}=\hat{L}^{ni}=0~~\text{for}~i<n.\eqno(3.7)$$

This change of variables is useful because, for any function $w$, $$\hat{L}^{ij}\p_{x_ix_j}(w(x',-x_n,t))=\hat{L}^{ij}(\p_{x_ix_j}w)(x',-x_n,t),$$
which will allow us to symmetrise the solution and still have a supersolution for the Pucci extremal operator $M^-$. We also use it to prove a semiconcavity result from semiconvexity in the following proof of Proposition 3.2.

Before proving Proposition 3.2, let us show the following lemma that is similar to Lemma 9.2 of \cite{C2}. The following lemma gives some properties of the linearized operator $L$ of $F$ at $D^2u.$
\begin{lemma}
Let $F$ be convex, $F\in C^\infty$ satisfy $F(0)=0$. Let $Q$ be a domain in $\rr^{n+1}$ and  $u\in H^4(Q)$ satisfy $F(D^2u)-\p_tu=0$ in $Q$. Consider the uniformly elliptic operator in $Q$ $$Lv=a_{ij}(x,t)v_{ij}:=F_{ij}(D^2u(x,t))v_{ij}.$$
Then, for any $e\in S^{n-1}$, $$Lu\geq \p_tu,\quad Lu_{e}=\p_tu_e\quad \text{and}\quad Lu_{ee}\leq \p_tu_{ee}~\text{in}~Q.$$
\end{lemma}
\begin{proof}
Let $\psi(s)=F((1-s)D^2u(x,t))-(1-s)\p_tu(x,t)$. We have that $\psi(0)=0$ and $\psi(1)=F(0)=0.$ By the convexity of $F$, it follows that $\psi\leq 0$ in $[0,1].$
Thus, $0\geq\psi'(0)=F_{ij}(D^2u(x,t))(-u_{ij}(x,t))+\p_tu(x,t)=-Lu(x,t)+\p_tu(x,t).$ That is, $Lu\geq \p_tu$ in $Q.$

We now differentiate $F(D^2u(x,t))-\p_tu(x,t)=0$ with respect to $e$ and get $Lu_e=\p_tu_e$. We differentiate it once more with respect to $e$, to obtain $$0=Lu_{ee}+F_{ij,kl}(D^2u(x,t))(u_e)_{ij}(u_e)_{kl}-\p_tu_{ee}.$$ By the convexity of $F$,  $Lu_{ee}\leq \p_tu_{ee}$ in $Q$. This completes the proof of the lemma.
\end{proof}
We next prove the semiconvexity of solutions in the directions parallel to the domains of the obstacle. To do it, we apply a Bernstein's technique in the spirit of \cite{AC}.

\begin{proposition}
Let $u$ be the solution to $(1.2)$. Then \\
\\
$(a)$~$(\mathrm{Semiconvexity})$ If $\tau=(\tau^*,0)$, with $\tau^*$ a unit vector in $\rr^{n-1},$
$$\inf_{Q_{3/4}} u_{\tau\tau}\geq -C\Big(\|u\|_{L^\infty(Q_1)}+\|\v\|_{H^{2}(Q^*_1)}\Big),$$
for some universal constant $C$.\\
$(b)$~$(\mathrm{Semiconcavity})$ Similarly, in the direction normal to $Q^*_1\times\{0\}$,
$$\sup_{Q_{3/4}}u_{x_nx_n}\leq C\Big(\|u\|_{L^\infty(Q_1)}+\|\v\|_{H^{2}(Q^*_1)}\Big),$$
for some universal constant $C$.
\end{proposition}
\begin{proof}
The second part, (b), follows from (a) by applying the definition of uniformly elliptic operator and the fact that we changed variables in order to have matrix $\hat{L}$ satisfying (3.7). We denote by $\hat{L}'$ and $D^2_{n-1}u$ the square matrices  corresponding to the $n-1$ first indices of $\hat{L}$ and $D^2u$ respectively. In fact, from $$\hat{L}^{ij}\p_{x_ix_j}u(x,t)\leq 0,\quad \hat{L}^{in}=\hat{L}^{ni}=0~~\text{for}~i<n,$$
and $$D^2_{n-1}u\geq -C\Big(\|u\|_{L^\infty(Q_1)}+\|\v\|_{H^{2}(Q^*_1)}\Big)\text{Id}_{n-1},$$
we can obtain that $$\hat{L}^{nn}\p_{x_nx_n}u\leq \p_tu-\sum_{i,j=1}^{n-1}\hat{L}^{ij}\p_{x_ix_j}u\leq \p_tu +C\Big(\|u\|_{L^\infty(Q_1)}+\|\v\|_{H^{2}(Q^*_1)}\Big)\text{tr}\hat{L}'.$$
The desired bound follows because $\hat{L}^{nn}$ is bounded below by $\l,$ $\mathrm{tr}(\hat{L}')$ is bounded above by $(n-1)\L$, and $\p_tu$ is locally upper bounded.

Now, we begin to prove (a). As the proof of Proposition 3.1, we define $h$ as the solution to (3.2) and (3.3). Recall that $h$ is Lipschitz and that, by maximum principle, $u>h$ in $Q^+_{1/2}$ and $Q^-_{1/2}.$

For $\varepsilon>0$, define $$\bar{h}_\varepsilon(x',x_n,t):=\v(x',t)-\frac{x^2_n}{\varepsilon}$$
and $$h_\varepsilon(x',x_n,t):=\max\{h(x',x_n,t),\bar{h}_\varepsilon(x',x_n,t)\}.$$
Since $h$ is Lipschitz continuous and $h(x',0,t)=\bar{h}_\varepsilon(x',0,t)$, this implies that there exists a universal constant $C>0$ such that $$h(x',x_n,t)>\bar{h}_\varepsilon(x',x_n,t)\quad\text{for}~|x_n|>CK\varepsilon.\eqno(3.8)$$
In particular, $h_\varepsilon$ is Lipschitz continuous in $Q_{7/8}$, uniformly on $\varepsilon.$

Let $u_\varepsilon$ be the solution to the parabolic ``thick'' obstacle problem with obstacle $h_\varepsilon$,
\begin{equation*}\label{1}
    \left\{
   \begin{array}{ll}
F(D^2u_\varepsilon)-\p_tu_\varepsilon=0\qquad &\text{in}~Q_1\setminus\{u_\varepsilon=h_\varepsilon\}\\
F(D^2u_\varepsilon)-\p_tu_{\varepsilon}\leq 0&\text{on}~Q_1\tag{3.9}\\
u_\varepsilon=\max\{u,\bar{h}_\varepsilon\}&\text{on}~ \p_pQ^+_1\setminus Q^*_1\\
u_\varepsilon\geq h_\varepsilon&\text{in}~Q^+_1,
 \end{array}
 \right.
\end{equation*}
and the analogous expression in $Q^-_1.$ It follows from (3.8) that the coincidence set satisfies $$\{u_\varepsilon=h_\varepsilon\}\subset\{\bar{h}_\varepsilon>h\}\subset\{(x',x_n,t)\in Q_1:|x_n|\leq CK\varepsilon\},$$
for some constant $C>0.$

Next, we want to bound $\p_{\tau\tau}u_\varepsilon$ from below independently of $\varepsilon.$ since $u_\varepsilon\geq h_\varepsilon,$ this also occurs along the free boundary, we obtain that $D^2(u_\varepsilon-h_{\varepsilon})\geq 0$ in the coincidence set. By the definition of $\bar{h}_\varepsilon$, this leads to $\p_{\tau\tau}u_\varepsilon\geq -CK$ in $\{u_\varepsilon=h_\varepsilon\}\cap Q_{7/8},$ for some universal constant $C>0.$ Therefore, it is enough to check that $\p_{\tau\tau}u_\varepsilon$ is uniformly bounded from below outside the coincidence set. We use a Bernstein's technique to proceed the proof.

Let $\eta\in C_c^\infty(Q_{7/8})$ be a smooth, cutoff function, with $0\leq \eta\leq 1$ and $\eta\equiv 1$ in $Q_{3/4}$. Define $$f_\varepsilon(x,t)=\eta(x,t)\p_{\tau\tau}u_\varepsilon(x,t)-\mu|\nabla_x u_\varepsilon(x,t)|^2$$ for some constant $\mu$ to be chosen later.
Since $h_\varepsilon$ is Lipschitz continuous independently of $\varepsilon$ in $Q_{7/8}$, then $|\nabla u_\varepsilon(x,t)|$ is bounded independently of $\varepsilon$ in $Q_{7/8}$. We split into three cases.\\
\\
$\mathrm{Case}~1:$ If the minimum $(x_0,t_0)$ in $Q_{7/8}$ is attained in the coincidence set, then $\p_{\tau\tau}u_\varepsilon(x_0,t_0)\geq -CK$ and we have that for every $(x,t)\in Q_{3/4}$,
\begin{align*}
\p_{\tau\tau}u_\varepsilon(x,t)&\geq \p_{\tau\tau}u_\varepsilon(x_0,t_0)-\mu|\nabla u_\varepsilon(x_0,t_0)|^2+\mu|\nabla u_\varepsilon(x,t)|^2\tag{3.10}\\
&\geq -CK-\mu\|\nabla u_\varepsilon\|^2_{L^\infty(Q_{7/8})}.
\end{align*}
$\mathrm{Case}~2:$ If the minimum $(x_0,t_0)$ is attained at the boundary, $\p_p Q_{7/8}$, then for every $(x,t)\in Q_{3/4},$
\begin{align*}
\p_{\tau\tau}u_\varepsilon(x,t)&\geq-\mu|\nabla u_\varepsilon(x_0,t_0)|^2+\mu|\nabla u_\varepsilon(x,t)|^2\tag{3.11}\\
&\geq -\mu\|\nabla u_\varepsilon\|^2_{L^\infty(Q_{7/8})}.
\end{align*}
$\mathrm{Case}~3:$ If the minimum $(x_0,t_0)$ of $f_\varepsilon$ in $Q_{7/8}$ is attained at some interior point outside the coincidence set $\{u_\varepsilon=h_\varepsilon\}$.

Let us assume that the operator $F$ is not only convex, but also $F\in C^\infty,$ so that  solutions are $H^{4}$ outside the coincidence set; see the end of the proof for the general case $F$ Lipschitz. In this case, the linearised operator of $F$ at $D^2u(x_0,t_0),$
$$Lv=a_{ij}v_{ij}:=F_{ij}(D^2u_\varepsilon(x_0,t_0))v_{ij},$$
is uniformly elliptic with ellipticity constants $\l$ and $\L.$ In addition, for any $e\in S^{n-1},$ $$Lu_\varepsilon(x_0,t_0)\geq \p_tu_\varepsilon,\quad L\p_{e}u_\varepsilon=\p_{t}\p_e{u_\varepsilon},\quad\text{and} \quad L\p_{ee}u_\varepsilon\leq \p_{t}\p_{ee}u_\varepsilon \quad \text{for}~(x_0,t_0).\eqno(3.12)$$

For convenience, in the following computations we denote $w=u_\varepsilon$. If $(x_0,t_0)$ is an interior minimum of $f_\varepsilon$ that is a $C^2$ function in $Q_{7/8}$, then $$0=\nabla_xf_\varepsilon(x_0,t_0)=(\nabla\eta w_{\tau\tau}+\eta\nabla w_{\tau\tau}-2\mu w_i\nabla w_i)(x_0,t_0),\eqno(3.13)$$
and $$0=\p_tf_\varepsilon(x_0,t_0)=(\eta_t w_{\tau\tau}+\eta w_{\tau\tau t}-2\mu w_i w_{it})(x_0,t_0).\eqno(3.14)$$
By (3.12), (3.13) and the fact that $(a_{ij})$ is elliptic, we deduce
\begin{align*}
0\leq a_{ij}f_{\varepsilon,ij}(x_0,t_0)&=L(f_\varepsilon)(x_0,t_0)\\
&\leq \big(a_{ij}\eta_{ij}w_{\tau\tau}+2a_{ij}\eta_i w_{\tau\tau,j}-\eta_t w_{\tau\tau}-2\mu a_{ij} w_{kj}w_{ki}\big)(x_0,t_0).\tag{3.15}
\end{align*}

Combining (3.13), (3.14) and (3.15), $$0\leq \Big(\Big(a_{ij}\eta_{ij}-\dfrac{2a_{ij}\eta_i\eta_j}{\eta}-\eta_t\Big)w_{\tau\tau}-2\mu a_{ij}w_{kj}w_{ki}+\frac{4\mu a_{ij}\eta_iw_{kj}w_k}{\eta}\Big)(x_0,t_0).\eqno(3.16)$$
Since $\sqrt{\eta}$ is Lipschitz, $|\nabla_{x,t} \eta|^2\leq C\eta.$ Thus, $$0\leq \big(C_0|w_{\tau\tau}|+\mu C_1|D^2w||\nabla w|-2\mu a_{ij}w_{kj}w_{ki}\big)(x_0,t_0),$$
for some constants $C_0$ and $C_1$ depending only on $n$ and $\L.$
From the uniform ellipticity of $(a_{ij}),$ it follows that  $$a_{ij}w_{ki}w_{kj}\geq \l C(n)|D^2w|^2.$$
Moreover, using $|w_{\tau\tau}(x_0,t_0)|\leq |D^2w(x_0,t_0)|$, we have $$|D^2w(x_0,t_0)|\leq \frac{C_0}{\mu}+C_1|\nabla w(x_0,t_0)|,$$
for some constants $C_0$ and $C_1$ depending only on $n, \l, \L.$
Since $(x_0,t_0)$ is a minimum in $Q_{7/8}$, for any $(x,t)\in Q_{3/4}$,
\begin{align*}
w_{\tau\tau}(x,t)&\geq \eta(x_0,t_0)w_{\tau\tau}(x_0,t_0)-\mu|\nabla u_\varepsilon(x,t)|^2\\
&\geq -|D^2w(x_0,t_0)|-\mu\|\nabla u_\varepsilon\|^2_{L^\infty(Q_{7/8})}\tag{3.17}\\
&\geq -\frac{C_0}{\mu}-C_1\|\nabla u_\varepsilon\|_{L^\infty(Q_{7/8})}-\mu\|\nabla u_\varepsilon\|^2_{L^\infty(Q_{7/8})}.
\end{align*}

We now choose $\mu=\|\nabla u_\varepsilon\|^{-1}_{L^\infty(Q_{7/8})}$. Notice that, in all three cases (3.10), (3.11) and (3.17),
we obtain $$\inf_{Q_{3/4}}\p_{\tau\tau}u_\varepsilon\geq -C\Big(\sup_{Q_{7/8}}|\nabla u_\varepsilon|+K\Big),$$
for some universal constant $C$. Since $u_\varepsilon$ is controlled by the Lipschitz norm of $u$ and Lipschitz continuous independently of $\varepsilon>0$, so that by Proposition 3.1, we conclude
\begin{align*}
\inf_{Q_{3/4}}\p_{\tau\tau}u_\varepsilon&\geq -C\Big(\|u\|_{\text{Lip}(Q_{7/8})}+\|\v\|_{H^{2}(Q^*_1)}+K\Big)\\
&\geq -C\Big(\|u\|_{L^\infty(Q_{1})}+\|\v\|_{H^{2}(Q^*_1)}\Big)\tag{3.18}.
\end{align*}

If $F$ is not smooth, then it can be regularised convoluting with a mollifier in the space of symmetric matrices, so that it can be approximated uniformly in compact sets by a sequence $\{F_k\}_{k\in \mathbb{N}}$ of convex smooth uniformly elliptic operators with ellipticity constants $\l$ and $\L.$ We can assume $F_k(0)=0$ by subtracting $F_k(0)$.
Since the obstacle $h_\varepsilon$ is $H^{2}$ in a neighbourhood of the free boundary, in $Q_{7/8}$ and for every $\varepsilon>0$ we have uniform $H^{1+\g}$ estimates in $k$ for solutions to (3.9) with operators $F_k$. Applying Arzel\`{a}-Ascoli theorem, that there exists a subsequence converging uniformly. Thus, the estimate (3.18) can be extended to solutions of (3.9) with operators not necessarily smooth. Therefore, (3.18) follows for any $F$ not necessarily $C^\infty.$

Note that $u_\varepsilon$ converges uniformly to $u$, since for all $\delta>0$, there exists some $\varepsilon>0$ small enough such that $u\leq u_\varepsilon<u+\delta$ in $Q_1.$

Since the right-hand side of (3.18) is independent of $\varepsilon$, and $u_\varepsilon$ converges uniformly to $u$ in $Q_{7/8}$ as $\varepsilon\downarrow 0$, we obtain $$\inf_{Q_{3/4}}u_{\tau\tau}\geq -C\Big(\|u\|_{L^\infty(Q_1)}+\|\v\|_{H^{2}(Q^*_1)}\Big).\eqno(3.19)$$  This completes the proof of the Lemma.
\end{proof}
\section{$H^{1+\alpha}$ estimate}
\subsection{A symmetrised solution}
By proposition 3.1, we know that $\nabla u$ is bounded in the interior of $Q_1$. In addition, $u_{x_nx_n}$ is bounded from Proposition 3.2 inside $Q_1$. In particular, for any $(x',t)\in Q^*_1$, the following limit exists $$\sigma (x',t)=\displaystyle{\lim_{x_n\downarrow 0^+}}u_{x_n}(x',x_n,t)-\displaystyle{\lim_{x_n \uparrow 0^-}}u_{x_n}(x',x_n,t)=\lim_{x_n\downarrow 0^+}\big(u_{x_n}(x',x_n,t)-u_{x_n}(x',-x_n,t)\big).\eqno(4.1)$$
An important step towards Theorem 1.1 consists of proving that $\s\in H^{\alpha}(Q^*_{1/2})$ for some $\alpha>0.$ We will prove it in this section.
\begin{lemma}
The function $\s$ defined by ${(4.1)}$ is nonpositive, i.e., $\s\leq 0$ in $Q^*_1.$
\end{lemma}
\begin{proof}
Let us argue by contradiction. Suppose that there exists some $(\bar{x}',\bar{t})\in Q^*_1$ such that $\s(\bar{x}',\bar{t})>0.$ Let $\d >0$ be such that $Q^*_{\d}(\bar{x}',\bar{t})\subset Q^*_1,$ so that by the semiconcavity in Proposition 3.2 applied to $Q_{\d/2}((\bar{x}',0,\bar{t}))$, $u_{x_nx_n}(\bar{x}',0,\bar{t})\leq C$ for some constant $C$, which now depends also on $\d$. However, $$\s(\bar{x}',\bar{t})=\lim_{x_n\downarrow 0^+}\big(u_{x_n}(\bar{x}',x_n,t)-u_{x_n}(\bar{x}',-x_n,t)\big)>0,$$
which implies $$\dfrac{u_{x_n}(\bar{x}',x_n,t)-u_{x_n}(\bar{x}',-x_n,t)}{2x_n}\rightarrow \iy,\quad \text{as}~x_n\downarrow 0^+.$$ This leads a contradiction with the bound in $u_{x_nx_n}.$
\end{proof}
We will now adapt the ideas of $\cite {C1}$ to our nonsymmetric setting. We define a symmetrised solution as follows $$v(x',x_n,t):=\dfrac{u(x',x_n,t)+u(x',-x_n,t)}{2}\quad \text{for}~ (x',x_n,t)\in \overline{Q_1}.\eqno(4.2)$$
Here $u$ is any solution to (1.2).

Notice that $$\s(x',t)=2\lim_{x_n\downarrow 0^+}v_{x_n}(x',x_n,t)\leq 0\eqno(4.3)$$
is well defined. Particularly, by the $H^{2+\alpha} $ interior estimates, we have that $$\s(x',t)=2v_{x_n}(x',0,t)=0 \quad \text{for} ~ (x',t)\in \Omega^*.\eqno(4.4)$$

The following result follows from the results in Section 3. We will denote $M^+$ and $M^-$ the Pucci extremal operators; see \cite{C2} for their definitions and basic properties.

\begin{lemma}
Let $u$ be a solution to the fully nonlinear parabolic thin obstacle problem $(1.2)$, $\p_t u$ is locally upper bounded, and let $v$ be defined by $(4.2)$. Then $v$ is spatial Lipschitz in $\overline{Q^+_{1/2}}$, $\p_t v$ is locally upper bounded, and satisfies
\begin{equation*}\label{1}
    \left\{
   \begin{array}{ll}
M^-(D^2v)-\partial_t v\leq 0 \qquad\qquad\quad&~\text{in}~Q_1,\tag{4.5}\\
\max\{v_{x_n}(x',0,t), \varphi(x',t)-v(x',0,t)\}=0&~\text{for}~(x',t)\in Q^*_1.
\end{array}\right.
\end{equation*}
Moreover,\\
\\
$(a)$~$(\mathrm{Semiconvexity})$ If $\tau=(\tau^*,0),$ with $\tau^*$ a unit vector in ${\rr}^{n-1},$
$$\inf_{Q_{3/4}} v_{\tau\tau}\geq -C\Big(\|u\|_{L^\infty(Q_1)}+\|\varphi\|_{H^{2+\alpha(Q^*_{1})}}\Big),$$
for some universal constant $C$.\\
$(b)$~$(\mathrm{Semiconcavity})$ In the direction normal to $Q^*_1\times\{0\}$,
$$\sup_{Q_{3/4}}v_{x_nx_n}\leq C\Big(\|u\|_{L^\infty(Q_1)}+\|\varphi\|_{H^{2}(Q^*_1)}\Big),$$
for some universal constant $C$.
\end{lemma}
\begin{proof}
The spatial Lipschitz regularity comes from the spatial Lipschitz regularity in $u$, proved in Proposition 3.1. $\p_t u$ is locally upper bounded. It is easy to see that $\p_t v$ is locally upper bounded.

Applying the change of variables introduced in section 3, we can obtain in (4.5) the first inequality. In fact, there exists some operator given by a matrix $\hat{L}$ as in (3.7) with uniformly elliptic constants $\l$ and $\L$ such that $$\hat{L}^{ij}\p_{x_ix_j}(u(x',-x_n,t))=\hat{L}^{ij}(\p_{x_ix_j}u)(x',-x_n,t)\leq F((D^2u(x',-x_n,t)).$$
Thus, $M^-(D^2v)-\p_tv\leq \hat{L}^{ij}\p_{x_ix_j}v-\p_tv\leq 0.$

The second equation in (4.5) follows from Lemma 4.1, (4.3), (4.4) and the fact that $v(x',0,t)=u(x',0,t)$ for $(x',t)\in Q^*_1.$

Finally, the semiconvexity and semiconcavity can be obtained by Proposition 3.2.
\end{proof}
\subsection{Regularity for $\s$ on free boundary points}
In this section, we follow the ideas of \cite {CG}, but we adapt them to the symmetrised solution $v$ instead of $u$. For the sake of completeness, we provide the detailed proofs. We begin with the following lemma, corresponding to Lemma 11 in \cite{CG}.

In the next result, we call $\varphi$ the extension of the obstacle to $Q_1$, i.e., $\varphi(x',x_n,t):=\varphi(x',t).$

\begin{lemma}
Let $v$ be the symmetrised solution $(4.2)$. Let $\k$ be a constant such that $\k > K$. For $(x_0,t_0)\in \Omega$ fixed, we define the function
$$\psi_{(x_0,t_0)}(x,t)=\v(x_0,t_0)+\nabla{\v}(x_0,t_0)\cdot (x-x_0)-\k(t-t_0)+\k|x-x_0|^2-\k n\frac{\lambda}{\Lambda}x_n^2.$$
We consider any set of the form $U:=U_{x_0}\times(t_1,t_0]\subset Q_1,$ with $U_{x_0}\subset \rr^n$ a bounded domain containing $x_0$ and $t_1<t_0<0$. Then
$$\sup_{\p_pU\cap\{x_n\geq 0\}}(v-\psi_{(x_0,t_0)})\geq 0.$$
\end{lemma}
\begin{proof}
Let $w:=v-\psi_{(x_0,t_0)}$. Using the definition of $\psi_{(x_0,t_0)}$ and Lemma 4.2, we obtain $w(x_0,t_0)\geq 0$ and $M^-(D^2w)-\p_t w\leq M^-(D^2v)+M^+(-D^2\psi_{(x_0,t_0)})-\p_tw\leq 2\k(\l-\L)-\k\leq 0$. Thus, we can apply the maximum principle on $U\setminus \bigtriangleup$ (recall $\bigtriangleup$ is the coincidence set) and use the symmetry of $w$ to obtain that $$\sup_{\p_p(U\setminus \triangle)\cap \{x_n\geq 0\}}(v-\psi_{(x_0,t_0)})\geq 0.$$

Now note that on the set $\{v=u=\v\}\cap\{t\leq t_0\}$ we have that $\psi_{(x_0,t_0)}> \v,$ since $(x_0,t_0)\in \Omega$, Remark 10 in \cite{CG} and $\k> \|\varphi\|_{H^{2}}.$ Therefore, $v-\psi_{(x_0,t_0)}<0$ on this set. Finally we observe that $\p_p(U\setminus\triangle)\cap\{x_n\geq 0\}\subset \p_pU\cap\{x_n\geq 0\}\bigcup(\triangle\cap \{t\leq t_0\}),$
then $$\sup_{\p_p(U\setminus \triangle)\cap \{x_n\geq 0\}}(v-\psi_{(x_0,t_0)})=\sup_{\p_pU\cap\{x_n\geq 0\}}(v-\psi_{(x_0,t_0)})\geq 0.$$ We complete the proof of this lemma.
\end{proof}

\begin{lemma}
Let $v$ be the symmetrised solution as defined in $(4.2),$ and let $\s$ as defined in $(4.1)$-$(4.3)$. For $\gamma>0$ we define $\Omega^*_{\gamma}:=\{(x',t)\in Q^*_1: \s(x',t)>-\gamma \}.$ Let $(x'_0,t_0)\in \Omega^*\cap Q^*_{1/2},$ then there exist constants $0<C<\bar{C}<1$ which depend only on $K, n, \l, \L$ so that for any $0<\gamma<\frac{1}{2}$ there exists a thin cylinder $Q^*_{C\gamma}(\bar{x}',\bar{t})$ such that $$Q^*_{C\gamma}(\bar{x}',\bar{t})\subset Q^*_{\bar{C}\gamma}(x'_0,t_0)\cap \Omega^*_{\gamma}.$$
\end{lemma}
\begin{proof}
Let $(x'_0,t_0)\in \Omega^*\cap Q^*_{1/2},$ we apply Lemma 4.3 with $$U:=B^*_{C_1\gamma}(x'_0)\times (-C_2\gamma,C_2\g)\times (t_0-(C_1\g)^2,t_0]$$
for some constants to be chosen $0<C_2<<C_1<<1$, and study two cases.

$\mathrm{Case}~1.$ Assume $\sup(v-\psi_{(x'_0,0,t_0)})$ is attained at a point $(x'_1,y_1,t_1)$ (for $x'_1\in \rr^{n-1}, y_1\in \rr $ and $t_1\in \rr$) on the lateral face of the cylinder $U$, i.e., with $|x'_0-x'_1|=C_1\g$ or $t_1=t_0-(C_1\g)^2$. Then applying Lemma 4.3 (for $\k=2K$) and Remark 10 in \cite{CG} we have
\begin{align*}
\psi_{(x'_0,0,t_0)}(x'_1,y_1,t_1)-\v(x'_1,t_1)&\geq K|x'_1-x'_0|^2-K(t_1-t_0)-2Kn\frac{\l}{\L}y^2_1\\
&\geq KC^2_1\g^2-2Kn\frac{\l}{\L}C^2_2\g^2\geq C_3\g^2,
\end{align*}
provided that $C_1>>C_2.$ The positive constant $C_3$ depends only on $K,n$, the ellipticity constants, $C_1$ and $C_2$. Therefore,
$$v(x'_1,y_1,t_1)\geq \psi_{(x'_0,0,t_0)}(x'_1,y_1,t_1)\geq \v(x'_1,t_1)+C_3\g^2.$$

Now pick a $(x'_2,t_2)\in Q^*_{C_4\g}(x'_1,t_1)$ for some positive $C_4$ to be chosen and $(x'_2-x'_1)\cdot\nabla_{x'}(v-\v)(x'_1,y_1,t_1)\geq 0$
(consider the extension of $\v$ in $Q_1$ where $\v(x',y,t)=\v(x',t)$). Take $\tau=\frac{x'_2-x'_1}{|x'_2-x'_1|}.$
We obtain
\begin{align*}
\int^{|x'_2-x'_1|}_0&\int^e_0(v-\v)_{\tau\tau}(x'_1+\tau h,y_1,t_1)dhde\\
=&(v-\v)(x'_2,y_1,t_1)-(v-\v)(x'_1,y_1,t_1)-|x'_2-x'_1|(v-\v)_\tau(x'_1,y_1,t_1),
\end{align*}
and $$\int^{t_1}_{t_2}(v-\v)_t(x'_2,y_1,h)dh=(v-\v)(x'_2,y_1,t_1)-(v-\v)(x'_2,y_1,t_2).$$
Combining the above we deduce
\begin{align*}
&\int^{|x'_2-x'_1|}_0\int^e_0(v-\v)_{\tau\tau}(x'_1+\tau h,y_1,t_1)dhde-\int^{t_1}_{t_2}(v-\v)_t(x'_2,y_1,h)dh\\
&=(v-\v)(x'_2,y_1,t_2)-(v-\v)(x'_1,y_1,t_1)-|x'_2-x'_1|(v-\v)_\tau(x'_1,y_1,t_1).\tag{4.6}
\end{align*}
Using $(a)$ of Lemma 4.2 we get
$$\int^{|x'_2-x'_1|}_0\int^e_0(v-\v)_{\tau\tau}(x'_1+\tau h,y_1,t_1)dhde\geq -C|x'_2-x'_1|^2\geq -C(C_4\g)^2$$
and $$-\int^{t_1}_{t_2}(v-\v)_t(x'_2,y_1,h)dh\geq -C(t_1-t_2)\geq -C(C_4\g)^2.$$
Thus returning to (4.6) we have
\begin{align*}
(v-\v)(x'_2,y_1,t_2)&\geq (v-\v)(x'_1,y_1,t_1)+(x'_2-x'_1)\cdot\nabla_{x'}(v-\v)(x'_1,y_1,t_1)-C(C_4\g)^2\\
&\geq C_3\g^2-C(C_4\g)^2>0\tag{4.7},
\end{align*}
if $0<C^2_4<\frac{C_3}{C}.$

To get a contradiction, now suppose that $(x'_2,t_2)\not\in \O^*_{\g},$ that is $\s(x'_2,t_2)\leq -\g<0.$ In particular, this implies $v(x'_2,0,t_2)=\v(x'_2,t_2).$
Similarly as before we want to transfer this information from $(x'_2,0,t_2)$ to $(x'_2,y_1,t_2)$ via integration of $v_{x_nx_n}$ and applying $(b)$ of Lemma 4.2, we obtain $$\int^{y_1}_0\int^e_0v_{x_nx_n}(x'_2,h,t_2)dhde=v(x'_2,y_1,t_2)-v(x'_2,0,t_2)-y_1\frac{\s(x'_2,t_2)}{2}\leq Cy^2_1.$$
Then, $v(x'_2,y_1,t_2)-\v(x'_2,t_2)\leq Cy^2_1-\frac{y_1\g}{2}\leq y_1\g(CC_2-\frac{1}{2})\leq 0$, provided that $0<C_2\leq \frac{1}{2C}.$
Thus, we have reached a contradiction regarding (4.7).\\
\\
$\mathrm{Case}~2.$ Assume now that $\sup (v-\psi_{(x'_0,0,t_0)})$ is attained at a point $(x'_1,y_1,t_1)$ in the base of the cylinder $U$, i.e., with $y_1=C_2\g$. Then applying Lemma 4.3 (for $\k=2K$) and Remark 10 in \cite{CG} we deduce
$$v(x'_1,y_1,t_1)\geq \psi_{(x'_0,0,t_0)}(x'_1,y_1,t_1)\geq \v(x'_1,t_1)-2Kn\frac{\l}{\L}C^2_2\g^2.$$
Now choose $(x'_2,t_2)\in Q^*_{C_2\g}(x'_1,t_1)$ and $(x'_2-x'_1)\cdot\nabla_{x'}(v-\v)(x'_1,y_1,t_1)\geq 0.$ We can modify slightly the computations of Case 1 to get $$(v-\v)(x'_2,C_2\g,t_2)\geq -(2Kn\frac{\l}{\L}+C)C^2_2\g^2>-C_5C_2\g^2,\eqno(4.8)$$
where $0<C_5<CC_2.$

To get a contradiction, now suppose that $(x'_2,t_2)\not\in \O^*_{\g},$ that is $\s(x'_2,t_2)\leq \g<0.$ Then $v(x'_2,0,t_2)=\v(x'_2,t_2)$. Similarly as Case 1 we deduce that $v(x'_2,C_2\g,t_2)-\v(x'_2,t_2)\leq C_2\g^2(CC_2-\frac{1}{2})<-C_5C_2\g^2$, provided that $0<C_5<\frac{1}{2}-CC_2$ and $C_2$ small enough. This is a contradiction regarding (4.8).

In any case, there exists $0<C_6<<1$ depending only on $K,n,\l$ and $\L$ so that if $(x'_2,t_2)\in Q^*_{C_6\g}(x'_1,t_1)$ with $(x'_2-x'_1)\cdot\nabla_{x'}(v-\v)(x'_1,y_1,t_1)\geq 0$ (which, roughly speaking, holds at least in the ``half'' of $Q^*_{C_6\g}(x'_1,t_1)$) then $(x'_2,t_2)\in \O^*_{\g}$. Furthermore, choose $1\geq \bar{C}>C_7+C_1$, and it is easy to check that $Q^*_{C_6\g}(x'_1,t_1)\subset Q^*_{\bar{C}\g}(x'_0,t_0)$. Finally, choosing a thin cylinder $Q^*_{{C}\g}(\bar{x}',\bar{t})$ inside $Q^*_{C_6\g}(x'_1,t_1)\cap\{(x'_2-x'_1)\cdot\nabla_{x'}(v-\v)(x'_1,y_1,t_1)\geq 0\}$, we complete the proof.
\end{proof}

In order to make an iterative argument to prove the $H^\alpha$ regularity of $\s$, we need the following lemma about solutions to the Pucci equations, and can be found in \cite{CG}, Lemma 13. It follows from an appropriate use of the strong maximum principle for $M^-$ and a barrier argument.

\begin{lemma}
Let the set $K_1:=B^*_1\times(0,1)\times(-1,0].$ Let $w$ be a nonnegative continuous function in $K_1$ that solves $$M^-(D^2w)-\p_tw\leq 0\quad \text{in}~K_1.$$
Suppose that there exists some neighborhood $Q^*_{\d}(\bar{x}',\bar{t})\subset Q^*_1$ so that $$\liminf_{x_n\downarrow 0^+} w(x',x_n,t)\geq 1 \quad \text{for any}~~(x',t)\in \overline{Q}^*_\d(\bar{x}',\bar{t}).$$
Then$$w(x',x_n,t)\geq \varepsilon>0\quad \text{for}~(x',x_n,t)\in \overline{B}^*_{1/2}\times\bigg[\frac{1}{4},\frac{3}{4}\bigg]\times\bigg[-\frac{\d^2}{2},0\bigg],$$
for some $\varepsilon$ depending only on $\d$, $n$, and the ellipticity constants $\l$ and $\L$.
\end{lemma}

We now show the following lemma, which is a consequence of an iterative argument, analogous to Lemma 14 of \cite{CG}.

\begin{lemma}
Let $\s$ as defined in $(4.1)$-$(4.3)$, for $u$ the solution to the parabolic thin obstacle problem $(1.2).$ Let $(x'_0,t_0)\in \O^*\cap Q^*_{1/2},$ then there exist universal constants $0<\alpha<1$ and $C>0$ such that $$\s(x',t)\geq -C\Big(\|u\|_{L^\infty(Q_1)+\|\v\|_{H^{2}(Q^*_1)}}\Big)\Big(|x'-x'_0|+|t-t_0|^{1/2}\Big)^\alpha\quad \text{for}~(x',t)\in Q^*_{1/2}(x'_0,t_0).$$

\end{lemma}
\begin{proof}
Define $$K:=\|u\|_{L^\infty(Q_1)}+\|\v\|_{H^{2}(Q^*_1)},$$
and notice that by taking $\frac{u}{K}$ instead of $u$ if necessary we can assume $$\|u\|_{L^\infty(Q_1)}+\|\v\|_{H^{2}(Q^*_1)}\leq 1.$$
In fact, if $K\geq 1$ then $$F_K(D^2u):=\frac{1}{K}F(D^2(Ku)),$$
is a convex elliptic operator with ellipticity constants $\l$ and $\L$, and $\frac{u}{K}$ is a solution to the nonlinear parabolic thin obstacle problem for the operator $F_K-\p_t$ with obstacle $\frac{\v}{K}$. In this case, $$\|\frac{u}{K}\|_{L^\infty(Q_1)}+\|\frac{\v}{K}\|_{H^{2}(Q^*_1)}=1.$$ Therefore, we can assume $K\leq 1.$

Our aim is to show by induction that there exists a universal constant $C>0$ for any $k\in \mathbb{N}$, such that $$v_{x_n}(x',x_n,t)\geq -C\t^k\quad \text{for}~(x',x_n,t)\in Q^*_{\g^k}(x'_0,t_0)\times\{x_n\in (0,\g^k)\},\eqno(4.9)$$
where $0<\g<<\t<1$ to be chosen. We proceed by induction. For $k=1$, recall that $\s(x',t)=2\displaystyle\lim_{x_n\downarrow 0^+}v_{x_n}(x',x_n,t),$ and that from Lemma 4.2, $v_{x_n}$ is bounded and $v_{x_nx_n}\leq C$. In addition, $\s$ is nonpositive by Lemma 4.1, so that $v_{x_n}\leq Cx_n,$ for $x_n>0.$ We assume that $(4.9)$ holds for some $k$ and we prove it for $k+1.$

Consider the function $$w:=\dfrac{v_{x_n}+C\t^k}{-\frac{\mu}{2} \g^k+ C\t^k}\quad \text{in} ~Q^*_{\g^k}(x'_0,t_0)\times\{x_n\in (0,\g^k)\},$$
where $0<\mu<\frac{1}{2}$ a small constant to be chosen. Using the hypothesis of the induction and choosing $\mu<2C$ and $\g<\t$, we get $w\geq 0$ in $Q^*_{\g^k}(x'_0,t_0)\times\{x_n\in (0,\g^k)\}.$ In addition, $M^-({D^2 w})-\p_tw\leq 0$ in $Q^*_{\g^k}(x'_0,t_0)\times\{x_n\in (0,\g^k)\}.$ For $(x',t)\in Q^*_{\g^k},$
$$\lim_{x_n\downarrow 0^+}w(x',x_n,t)=\dfrac{\frac{\s(x',t)}{2}+C\t^k}{-\frac{\mu}{2}\g^k+C\t^k}.$$
Applying Lemma 4.4 around $(x'_0,t_0)\in \O^*\cap Q^*_{1/2}$, we obtain that there exists $Q^*_{\bar{C}\mu\g^k}(\bar{x}',\bar{t})\subset Q^*_{\mu\g^k}(x'_0,t_0)\cap \O^*_{\mu\g^k}$, where $0<\bar{C}<1$ depends only on $K,n,\l$ and $\L.$ That is, $$\lim_{x_n\downarrow 0^+}w(x',x_n,t)\geq 1 \quad \text{for}~ (x',t)\in Q^*_{\bar{C}\mu\g^k}(\bar{x}',\bar{t}).$$ Thus, $w$ fulfils the hypothesis of Lemma 4.5 in $Q^*_{\g^k}(x'_0,t_0)\times\{x_n\in (0,\g^k)\}.$ Now we apply Lemma 4.5 to the rescaled function $W(x',x_n,t):=w(\mu\g^kx'+x'_0,\mu\g^kx_n,(\mu\g^k)^2t+t_0)$ in $K_1$ and get that $$w\geq \varepsilon\quad \text{in}~\overline{B}^*_{\frac{\mu\g^k}{2}}(x'_0)\times\bigg[\frac{\mu\g^k}{4},\frac{3\mu\g^k}{4}\bigg]\times\bigg[t_0-\frac{(\bar{C}\mu\g^k)^2}{2},t_0\bigg],\eqno(4.10)$$
where $\varepsilon=\varepsilon(\bar{C},n,\l,\L)>0.$ Since $r<\t$ and we choose $\mu$ so that $\mu\leq C$, $$v_{x_n}\geq -C\t^k+\frac{\varepsilon C\t^k}{2}.$$

Now, by means of Lemma 4.2, $v_{x_nx_n}\leq C',$ we deduce that
\begin{align*}
v_{x_n}(x',x_n,t)&= -\int^{\frac{\mu\g^k}{2}}_{x_n}v_{x_nx_n}(x',h,t)dh+v_{x_n}(x',\frac{\mu\g^k}{2},t)\\
&\geq-C'(\frac{\mu\g^k}{2}-x_n)-C\t^k+\frac{\varepsilon C\t^k}{2}
\end{align*}
for any $(x',x_n,t)\in \overline{B}^*_{\frac{\mu\g^k}{2}}(x'_0)\times\Big(0,\frac{\mu\g^k}{4}\Big]\times\Big[t_0-\frac{(\bar{C}\mu\g^k)^2}{2},t_0\Big].$
Furthermore, $$v_{x_n}(x',x_n,t)\geq -C\t^k+\frac{\varepsilon C\t^k}{2}-\frac{C'\mu\g^k}{2}$$
for $(x',x_n,t)\in \overline{B}^*_{\frac{\mu\g^k}{2}}(x'_0)\times\Big(0,\frac{3\mu\g^k}{4}\Big]\times\Big[t_0-\frac{(\bar{C}\mu\g^k)^2}{2},t_0\Big].$
Choosing $0<\g<\min\{\frac{\mu}{2},\frac{\sqrt{2}\bar{C}\mu}{2}\}<\frac{1}{2}$, we have the above holds for $(x',x_n,t)\in \overline{B}^*_{\g^{k+1}}(x'_0)\times(0,\g^{k+1})\times[t_0-(\g^{k+1})^2,t_0].$ To end the inductive argument we must see $$C\t^{k+1}\geq C\t^k-\frac{\varepsilon C\t^k}{2}+\frac{C'\mu\g^k}{2}.$$
For this, we choose $\mu<\frac{C\varepsilon}{4C'}$ and $\t>1-\frac{\varepsilon}{4}.$ Thus, the inductive argument is completed.

Taking $x_n\rightarrow 0$ in $(4.10)$ we obtain that  for any $k\in \mathbb{N}$ $$\s(x',x_n,t)\geq -C\t^k\quad \text{for} ~(x',t)\in Q^*_{\g^k}(x'_0,t_0),$$
where $0<\g<<\t<1$ and $C=C(K,n,\l,\L)>0$. The desired regularity for $\s$ follows in a standard way.
\end{proof}

\subsection{Proof of Theorem 1.1}
Before proving Theorem 1.1, we show the following compactness lemma.

\begin{lemma}
Let $F$ be a fully nonlinear operator satisfying $(1.3)$, and let $w$ be a continuous function defined in $Q_1$. Assume that $w$ satisfies the problem $$F(D^2w)-\p_tw=0\quad \text{in}~Q^+_1\cup Q^-_1,\eqno(4.11)$$
and that $$\|w\|_{L^\infty(Q_1)}=1,\quad [w]_{\mathrm{Lip}(Q_1)}\leq 1.$$

Let $\psi$ be the solution to
\begin{equation*}\label{1}
    \left\{
   \begin{array}{ll}
F(D^2\psi)-\p_t\psi=0\quad &\text{in}~Q_1\tag{4.12}\\
\psi=w&\text{on}~\p_pQ_1,
\end{array}\right.
\end{equation*}
and let us define the following operator $$\bar{\s}(w):=\lim_{h_n\downarrow 0}\big((\p_{x_n}w)(x',h_n,t)-(\p_{x_n}w)(x',-h_n,t)\big).$$

Then, for every $\varepsilon>0$ there exists some $\eta=\eta(\varepsilon,n,\l,\L)>0$ such that if $$\|\bar{\s}(w)\|_{L^\infty(Q^*_1)}<\eta$$
then $$\|\psi-w\|_{L^\infty(Q_1)}<\varepsilon,$$
i.e., $\psi$ approximates $w$ as $\eta$ tends to 0.
\end{lemma}
\begin{proof}
Assume by contradiction that there exists some fixed $\varepsilon>0$, we have a sequence $\eta_k\rightarrow 0$, fully nonlinear convex operators uniformly elliptic with ellipticity constants $\l$ and $\L$ with $F_k(0)=0$ and functions $w_k$, such that $$F_{k}(D^2w_k)-\p_tw_k=0\quad \text{in}~Q^+_1\cup Q^-_1\eqno(4.13)$$
and $$\|w_k\|_{L^\infty(Q_1)}=1,\qquad[w_k]_{\mathrm{Lip}(Q_1)}\leq 1,$$ with $$\|\bar{\s}(w_k)\|_{L^\infty(Q^*_1)}<\eta_k, \eqno(4.14)$$
but such that $$\|\psi_k-w_k\|_{L^\infty(Q_1)}\geq \varepsilon\quad\text{for all}~k,\eqno(4.15)$$
where $\psi_k$ is the solution to
\begin{equation*}\label{1}
    \left\{
   \begin{array}{ll}
F_k(D^2\psi_k)-\p_t\psi_k=0\quad &\text{in}~Q_1\tag{4.16}\\
\psi_k=w_k&\text{on}~\p_pQ_1.
\end{array}\right.
\end{equation*}

Applying Arzel\`{a}-Ascoli, up to a subsequence, we obtain that $w_k$ converges to some function $\bar{w}$ uniformly in $Q_1$, with $\|\bar{w}\|_{L^\infty(Q_1)}=1.$
On the other hand, up to a subsequence, the operators $F_k$ converge to some fully nonlinear convex uniformly elliptic operator $\bar{F}$ with ellipticity constants $\l$ and $\L$ such that $\bar{F}(0)=0$. Note also that $\psi_k$ converges uniformly to the solution $\bar{\psi}$ to
\begin{equation*}\label{1}
    \left\{
   \begin{array}{ll}
\bar{F}(D^2\bar{\psi})-\p_t\bar{\psi}=0\quad &\text{in}~Q_1\tag{4.17}\\
\bar{\psi}=\bar{w}&\text{on}~\p_pQ_1.
\end{array}\right.
\end{equation*}
By (4.15), in the limit $$\|\bar{\psi}-\bar{w}\|_{L^\infty(Q_1)}\geq \varepsilon>0. \eqno(4.18)$$

Consider the function $w_k+\eta_k|x_n|$ in $Q_1$. By (4.14), we have $$\bar{\s}(w_k+\eta_k|x_n|)\geq \eta_k>0\quad \text{in}~Q^*_1.$$

Since $F_k(D^2w_k)-\p_tw_k=0$ in $Q^+_1\cup Q^-_1,$ we obtain that, in the viscosity sense, $$F_k(D^2(w_k+\eta_k|x_n|))-\p_t(w_k+\eta_k|x_n|)\geq 0\quad \text{in}~Q_1.$$

Now, passing to the limit, noticing that $w_k+\eta_k|x_n|$ converges uniformly to $\bar{w}$, we reach that, in the viscosity sense, $$\bar{F}(D^2\bar{w})-\p_t\bar{w}\geq 0\quad\text{in}~Q_1.$$

Next, consider the function $w_k-\eta_k|x_n|$. Repeating the argument for $w_k-\eta_k|x_n|$, we get $$\bar{F}(D^2\bar{w})-\p_t\bar{w}\leq 0\quad\text{in}~Q_1.$$

Combining the above we conclude $$\bar{F}(D^2\bar{w})-\p_t\bar{w}=0\quad\text{in}~Q_1.$$
This implies $\bar{\psi}=\bar{w}$ in $Q_1$. This is a contradiction regarding (4.18).
\end{proof}

We now give the proof of Theorem 1.1.\\
\\
$Proof~of~Theorem~1.1.$\quad We separate the proof into three steps. Firstly, we prove that the solution $u$ is $H^{1+\alpha}$ around points in $\O^*$ by means of Lemma 4.6 and Lemma 4.7. Secondly, we apply the result from the first step to deduce that $\s$ is $H^{\alpha}$ in $Q^*_{1/4}$. Finally, the proof is completed in the third step.

As in the proof of Lemma 4.6 we assume that $$\|u\|_{L^\infty(Q_1)}+\|\v\|_{H^{2}}(Q^*_1)\leq 1,$$
to avoid having this constant on each estimate throughout the proof.\\
\\
{\bf{Step 1:}} Assume that the origin is a free boundary point. Under this assumption we will prove that there exists some affine function $L=A+B\cdot x$ such that $$\|u-L\|_{L^\infty(Q_r)}\leq Cr^{1+\alpha},\quad \text{for~all}~r>0,\eqno(4.19)$$
for some universal constants $C$ and $\alpha>0$. We proceed in the spirit of the proof of Lemma 1.5 in \cite{W2}.

Note that from Lemma 4.1 and Lemma 4.6 we know that there exists $\eta>0$ such that  $$|\s(x',t)|\leq \eta(|x'|+|t|^{1/2})^\alpha\quad \text{for~all}~(x',t)\in Q^*_1.\eqno(4.20)$$
Up to replacing from beginning $u(x,t)$ by $u(r_0x,r_0^2t)$ with $r_0<<1$, we can make $\eta$ as small as necessary.

Now let us show that there exists $\rho=\rho(\alpha,n,\l,\L)<1$ and a sequence of affine functions $$L_k(x)=A_k+B_k\cdot x\eqno(4.21)$$
such that $$\|u-L_k\|_{L^\infty(Q_{\rho^k})}\leq \rho^{k(1+\alpha)},\eqno(4.22)$$
and $$|A_k-A_{k-1}|+\rho^{k}|B_k-B_{k-1}|\leq C\rho^{(k+1)(1+\alpha)},\eqno(4.23)$$
for some universal constant $C$.

We prove it by induction on $k$. The conclusion for $k=0$ is satisfied by taking $L_0=0$. Suppose that $k$-th step is true. Let $$w_k(x,t)=\dfrac{(u-L_k)(\rho^kx,\rho^{2k }t)}{\rho^{k(1+\alpha)}}\quad \text{for}~(x,t)\in Q_1.$$

Notice that $$F_k(D^2w_k)-\p_tw_k=0\quad \text{in}~Q^+_1\cup Q^-_1$$
for some operator $F_k$ of the form (1.3). From the induction assumption, $$\|w_k\|_{Q_1}\leq 1.$$

In addition, let $$\s_k(x',t)=\lim_{h \downarrow 0}\big(\p_{x_n}w_k(x',h,t)-\p_{x_n}w_k(x',-h,t)\big)\quad\text{for}~(x',t)\in Q^*_1.$$
Then, one can check that, by (4.20), $$|\s_k(x',t)|\leq \eta(|x'|+|t|^{1/2})^\alpha.$$

Given $\varepsilon>0$ small, by Lemma 4.7, we can choose $\eta$ small enough such that $$\|\psi_k-w_k\|_{L^\infty(Q_1)}\leq \varepsilon,$$
where $\psi_k$ is the solution to
\begin{equation*}\label{1}
    \left\{
   \begin{array}{ll}
F_k(D^2\psi_k)-\p_t\psi_k=0\quad &\text{in}~Q_1\tag{4.24}\\
\psi_k=w_k&\text{on}~\p_pQ_1.
\end{array}\right.
\end{equation*}

By interior estimate, $\psi_k$ is $H^{2+\beta}$ in $Q_{1/2}$ with estimates depending only on $n$, $\l$ and $\L.$ Then, let $l_k$
be the linearisation of $\psi_k$ around 0, so that up to choosing $\rho$,
\begin{align*}
\|w_k-l_k\|_{L^\infty(Q_\rho)}&\leq \|w_k-\psi_k\|_{L^\infty(Q_\rho)}+\|\psi_k-l_k\|_{L^\infty(Q_\rho)}\\
&\leq \varepsilon+C\rho^{2},
\end{align*}
where $l_k=\nabla\psi_k(0,0)\cdot x+\psi_k(0,0).$ Choosing $\rho$ so that  $C\rho^{2}\leq\frac{1}{2}\rho^{1+\alpha}$ and choosing $\eta$ so that $\varepsilon\leq \frac{1}{2}\rho^{1+\alpha}$, we obtain $$\|w_k-l_k\|_{L^\infty(Q_\rho)}\leq \rho^{1+\alpha}.$$

Now, recalling the definition of $w_k$, $$\bigg\|u-L_k-\rho^{k(1+\alpha)}l_k\bigg(\frac{\cdot}{\rho^k}\bigg)\bigg\|_{L^\infty(Q_{\rho^{k+1}})}\leq \rho^{(k+1)(1+\alpha)},$$
so that the inductive step is concluded by taking $$L_{k+1}(x)=L_k(x)+\rho^{k(1+\alpha)}l_k\bigg(\dfrac{x}{\rho^k}\bigg).$$
It is clear that $L_{k+1}$ satisfies the required conditions.

Once one has (4.21), (4.22) and (4.23), define $L$ as the limit of $L_k$ as $k\rightarrow \infty$ (which exists, by (4.23)). And notice that, given any $0<r=\rho^k$ for some $k\in\mathbb{N}$, then $$\|u-L\|_{L^\infty(Q_r)}\leq \|u-L_k\|_{L^\infty(Q_r)}+\displaystyle\sum_{i\geq k}\|L_{i+1}-L_i\|_{L^\infty(Q_r)}\leq Cr^{1+\alpha},$$
for some universal constant C.\\
\\
{\bf{Step 2:}} In this step we prove that the function $\s$ defined in (4.1)-(4.3) is $H^\alpha(Q^*_{1/4})$ for some $\alpha=\alpha(n,\l,\L)>0$, and $$\|\s\|_{H^\alpha(Q^*_{1/4})}\leq C,\eqno(4.25)$$
for some universal constant $C$.

We already know $\s$ is regular in the interior of $\triangle^*$ (by boundary estimates) and $\O^*$. In particular, from the interior estimates $\s=0$ in $\O^*$. Applying Lemma 4.6, we also obtain $H^\alpha$ regularity at points in $\p\triangle^*.$ That is, we have that given $(x'_0,0,t)\in\p\triangle^*,$
$$|\s(x',t)|\leq C(|x'-x'_0|+|t-t_0|^{1/2})^\alpha \quad \text{for}~(x',t)\in Q^*_{1/2},\eqno(4.26)$$
for some universal constant $C$.

Thus, we only need to check that given $(x,t), (y,s)\in \triangle$, $(x,t)=(x',0,t)$, $(y,s)=(y,0,s)$, then there exists some universal constant $C$ such that, if $p((x,t),(y,s))=r,$$$|\s(x',t)-\s(y',s)|\leq Cr^\alpha.$$

Let $R:=$dist$_p((x,t),\O)$ and suppose that dist$_p((x,t),\O)\leq$ dist$_p((y,s),\O).$ Let $(z,q)=(z',0,q)$, $(z',q)\in \p\triangle^*,$ be such that dist$_p((x,t),(z,q))=$dist${_p}((x,t),\O),$ and assume that $\displaystyle\lim_{x_n\downarrow 0^+}\nabla u(z',x_n,q)=0$ and $\nabla_{x'}\v(z',q)=0$ by subtracting an affine function if necessary. Note that we can do so because we already know from the first step that $u$ is $H^{1+\alpha}$ around $(z',q)$. We separate into two cases.\\
\\
$\mathrm{Case}~ 1:$ If $R<4r,$ then applying (4.26)
\begin{align*}
|\s(x',t)-\s(y',s)|&\leq |\s(x',t)-\s(z',q)|+|\s(z',q)-\s(y',s)|\\
&\leq C[R^\alpha+(R+r)^\alpha]\\
&\leq Cr^\alpha.
\end{align*}
$\mathrm{Case} ~2:$ In the case $R\geq 4r$. Notice that $(x',t), (y',s)\in Q^*_{R/2}(x',t)\subset Q^*_R(x',t)\subset \triangle^*,$ and $u$ restricted to $Q^*_R(x',t)$ is a $H^{2}$ function, since $u=\v$. In particular, we use that under these assumption $$R^{1+\alpha}[\nabla u]_{{\alpha};\overline{Q^+_{R/2}(x,t)}}\leq C\Big(\mathop{\text{osc}}\limits_{Q^+_R(x,t)}u+R^2[D^2\v]_{Q^*_R(x',t)}\Big);$$
see Theorem 17 in \cite{CG1}. Recall that the gradient of $u$ at $(z,q)$ is 0. Consequently, from the previous step using the bound (4.19) around $(z,q)$,
$$|u(y,s)-\v(z',q)|\leq C(|y-z|+|s-q|^{1/2})^{1+\alpha}\leq CR^{1+\alpha}\quad \text{for}~(y,s)\in Q^+_R(x,t).$$
Particularly, $\mathop{\text{osc}}\limits_{Q^+_R(x,t)}u\leq CR^{1+\alpha}.$ Thus, this leads $$[\nabla u]_{{\alpha};\overline{Q^+_{R/2}(x,t)}}\leq C,$$
from which (4.25) is proved.
\\
\\
{\bf{Step 3:}} Note that in the first step we only used that the origin was a free boundary point to be able to apply   (4.20) in Lemma 4.6 . Now, given any point $(z',q)\in Q^*_{1/4},$  we can consider the function $u_{(z,q)}(x,t):=u(x,t)-\s(z',q)(x_n)^+,$
where $(x_n)^+$ denotes the positive part of $x_n.$

Notice that this function satisfies the assumption of Step 1, for $(x',t)\in Q^*_1,$
$$|\s_{(z,q)}(x',t)|:=\Big|\lim_{h\downarrow 0}\big(\p_{x_n}u_{(z,q)}(x',h,t)-\p_{x_n}u_{(z,q)}(x',-h,t)\big)\Big|\leq C(|x'-z'|+|t-q|^{1/2})^\alpha,$$
for some universal constant $C.$

Repeating the same procedure as in Step 1, we obtain that for every point $(z,q)\in Q_{1/4}\cup\{x_n=0\}$, and for every $(x,t)\in Q^+_1$ there exists some $L^+_{(z,q)}$ affine function such that $$|u(x,t)-L^+_{(z,q)}|\leq C(|x-z|+|t-q|^{1/2})^{1+\alpha}.$$ On the other hand, for every $(x,t)\in Q^-_1$ there exists some $L^-_{(z,q)}$ affine function such that $$|u(x,t)-L^-_{(z,q)}|\leq C(|x-z|+|t-q|^{1/2})^{1+\alpha}.$$ Combining the above, we conclude $$\|u\|_{H^{1+\alpha}(Q^*_{1/4})}\leq C,$$
for some universal constant $C$.

To finish the proof, we could now repeat a procedure like the one done in Step 2. \qed

LMAM, School of Mathematical Sciences, Peking University, Beijing, 100871, P. R. China

Xi Hu,\quad
E-mail address: huximath1994@163.com

Lin Tang,\quad
E-mail address: tanglin@math.pku.edu.cn
\end{document}